\documentclass[12pt]{article}
\usepackage{amssymb}
\usepackage{latexsym,bm}
\usepackage{graphicx}
\usepackage{amsmath}
\usepackage{mathrsfs}

\newcommand{\bea}{\begin{eqnarray*}}
\newcommand{\eea}{\end{eqnarray*}}
\newcommand{\be}{\begin{equation}}
\newcommand{\ee}{\end{equation}}
\newcommand{\ben}{\begin{eqnarray*}}
\newcommand{\een}{\end{eqnarray*}}

\date{}
\begin{document}
\title{The maximum number of cliques in graphs with prescribed order, circumference and minimum degree\footnote{E-mail addresses:
{\tt mathdzhang@163.com}.}}
\author{\hskip -10mm Leilei Zhang\\
{\hskip -10mm \small Department of Mathematics, East China Normal University, Shanghai 200241, China}}\maketitle
\begin{abstract}
 Erd\H{o}s determined the maximum size of a nonhamiltonian graph of order $n$ and minimum degree at least $k$ in 1962. Recently, Ning and Peng generalized
 Erd\H{o}s' work and gave the maximum number of $s$-cliques $h_s(n,c,k)$ of graphs with prescribed order $n$, circumference $c$ and minimum degree at least $k.$ But for some quadruple $n,c,k,s$ the maximum number of $s$-cliques is not attained by a graph of minimum degree $k.$ For example, $h_2(15,14,3)=77$ is attained by a unique graph of minimum degree $7,$ not $3.$ In this paper we obtain more precise information by determining the maximum number of $s$-cliques of a graph with prescribed order, circumference and minimum degree. Consequently we solve the corresponding problem for longest paths.
\end{abstract}



\section{ Introduction}
We consider finite simple graphs, and use standard terminology and notations. The {\it order} of a graph is its number of vertices, and the
{\it size} is its number of edges. For graphs we will use equality up to isomorphism, so $G_1=G_2$ means that $G_1$ and $G_2$ are isomorphic.
$\overline{G}$ denotes the complement of a graph $G.$ For two graphs $G$ and $H,$ $G\vee H$ denotes the {\it join} of $G$ and $H,$ which is obtained from the disjoint union $G+H$ by adding edges joining every vertex of $G$ to every vertex of $H.$ $K_n$ denotes the complete graph of order $n.$

Denote by $V(G)$ and $E(G)$ the vertex set and edge set of a graph $G.$  Let $e(G)$ denote the size of $G.$ For a vertex $v$ in a graph, we denote by $d(v)$ and $N(v)$ the degree of $v$ and the neighborhood of $v$ in $G,$ respectively. For $H\subseteq V(G),$ we denote by $N_H(v)$ the set $H\cap N(v)$ and denote
$d_H(v)=|N_H(v)|.$ We denote by $\delta(G)$ the {\it minimum degree} of a graph $G.$  For two vertices $u$ and $v$, we use the symbol $u\leftrightarrow v$ to mean that $u$ and $v$ are adjacent and use $u\nleftrightarrow v$ to mean that $u$ and $v$ are non-adjacent.

In 1961 Ore [14] determined the maximum size of a nonhamiltonian graph with a given order and also determined the extremal graphs.

{\bf Theorem 1.} (Ore [14]) {\it The maximum size of a nonhamiltonian graph of order $n$ is $\binom{n-1}{2}+1$ and this size is attained by a graph $G$
if and only if $G=K_1\vee (K_{n-2}+K_1)$ or $G=K_2\vee \overline{K_3}.$}

Bondy [1] gave a new proof of Theorem 1. It is natural to ask the same question by putting constraints on the graphs. In 1962 Erd\H{o}s [5] determined
the maximum size of a nonhamiltonian graph of order $n$ and minimum degree at least $k.$

{\bf Theorem 2.} (Erd\H{o}s [5]) {\it Let $n$, $k$ be integers with $1\le k\le \lfloor \frac{n-1}{2}\rfloor,$ and set $h(n,k)=\binom{n-k}{2}+k^2.$ If $G$ is a nonhamiltonian graph of order $n$ with $\delta(G)\ge k$, then
$$
e(G)\le {\rm max}\{h(n,k),\, h(n,\left\lfloor\frac{n-1}{2}\right\rfloor)\}.
$$ }

F\"{u}redi, Kostochka and Luo [8] proved a stability version of this theorem in 2018. The {\it circumference} $c(G)$ of a graph $G$ is the length of a longest cycle in $G.$ Determining the circumference of a graph is a classical  problem in graph theory. It is well known that even determining if a graph is hamiltonian is NP-hard. One cornerstone in this direction is the following celebrated Erd\H{o}s-Gallai
theorem.

{\bf Theorem 3.} (Erd\H{o}s and Gallai [6]) {\it Let $G$ be a graph with order $n$ and circumference $c$. Then $e(G)\le c(n-1)/2.$}

This is sharp if $n-1$ is divisible by $c-1.$ This can be seen by considering the graph $(\frac{n-1}{c-1}K_{c-1})\vee K_1.$ Theorem 3 also implies that if an $n$-vertex graph $G$ contains no path of length $p$, then $e(G)\le (p-1)n/2.$  F\"{u}redi, Kostochka and Verstra\"{e}te [9] proved a stability version of Theorem 3.

{\bf Notation 1.} An {\it s-clique} is a clique of cardinality $s.$ Fix $n-1 \geq c\geq 2k\ge 4.$ Let $F(n,c,k)=K_k\vee(K_{c+1-2k}+\overline{K_{n-c-1+k}})$.
For $t=\lfloor c/2\rfloor$, let $G(n,c,k)$ denote the graph obtained from $K_t\vee(K_{c+1-2t}+\overline{K_{n-c-1+t}})$  by deleting $t-k$ edges that are incident to one common vertex in $\overline{K_{n-c-1+t}}$. Denote by $f_s(n,c,k)$ the number of $s$-cliques in $F(n,c,k)$ and by $g_s(n,c,k)$ the number of $s$-cliques in $G(n,c,k)$; more precisely,
\begin{align*}
  f_s(n,c,k)&=\binom{c+1-k}{s}+(n-c-1+k)\binom{k}{s-1},\\
  g_s(n,c,k)&=\binom{c+1-t}{s}+(n-c-2+t)\binom{t}{s-1}+\binom{k}{s-1}.
\end{align*}

We write $f(n,c,k)$ for $f_2(n,c,k)$ which equals the size of $F(n,c,k),$ and we write $g(n,c,k)$ for $g_2(n,c,k)$ which equals the size of $G(n,c,k).$

By imposing minimum degree as a new parameter, Woodall [16] asked the following refinement of Theorem 3 in 1976.

{\bf Conjecture 4.} (Woodall [16]) {\it Let $c\le n-1$. The maximum size of a $2$-connected graph with order $n,$ circumference $c$ and  minimum degree at least $k$ is
$$
{\rm max}\{f(n,c,k), \, f(n,c,\lfloor c/2 \rfloor)\}.
$$ }

Woodall's conjecture 4 has been proved (see [7],[11],[13],[16]). For further developments on this topic, see [12].
 Improving Theorem 3, Kopylov [10] proved the following result in 1977.

{\bf Theorem 5.} (Kopylov [10]) {\it The maximum size of a $2$-connected graph of order $n$ and circumference $c$ with $c\le n-1$ is
$$
{\rm max}\{f(n,c,2),\, f(n,c,\lfloor c/2\rfloor)\}.
$$}
Hence when $f(n,c,\lfloor c/2\rfloor)> f(n,c,2),$ the graph which attains this bound has minimum degree $\lfloor c/2\rfloor.$

{\bf Notation 2.} For $s\geq 2$, let $N_s(G)$ denote the number of $s$-cliques in $G$; e.g., $N_2(G)=e(G)$.

Generalizing Theorem 5 (but using Kopylov's proof idea), Luo [11] proved the following result.

{\bf Theorem 6.} (Luo [11]) {\it Let $n-1\geq c\geq 4$ and let $s\ge 2$. If $G$ is a $2$-connected $n$-vertex graph with $c(G)=c$, then
$$
N_s(G)\leq {\rm max}\{f_s(n,c,2),\, f_s(n,c,\lfloor c/2\rfloor)\}.
$$}

Ning and Peng present an extension of Theorem 6 by adding a restraint on the minimum degree.

{\bf Theorem 7.} (Ning and Peng [13]) {\it Let $c\le n-1$ and $s\geq 2$. If $G$ is a $2$-connected graph of order $n$ with  $c(G)=c$ and minimum degree $\delta(G)\geq k\geq 2$, then
$$
N_s(G)\leq {\rm max}\{f_s(n,c,k),\, f_s(n,c,\lfloor c/2\rfloor)\}.
$$}

Let $h_s(n,c,k)$ denote the maximum number of $s$-cliques in graphs with prescribed order $n$, circumference $c$ and minimum degree at least $k.$ One can easily find that Theorem 2 and Theorem 6 are special cases of theorem 7. But for some cases, the bound can not be attained by a graph of minimum degree $k.$ For example, $h_2(15,14,3)=77$ is attained  by a unique graph of minimum degree $7,$ not $3.$

In this paper we obtain more precise information by determining the maximum number of $s$-cliques of a graph with given order, circumference and minimum degree. Consequently we solve the corresponding problem for longest paths. One of our main results is the following:

{\bf Theorem 8.}  {\it Let $\varphi_s(n,c,k)$ denote the maximum number of $s$-cliques in a $2$-connected nonhamiltonian graph of order $n$ with circumference $c$ and minimum degree $k.$
 Then
$$
\varphi_s(n,c,k)= {\rm max}\{f_s(n,c,k),g_s(n,c,k)\}.
$$ }

In Section 2 we give a proof of Theorem 8, and in Section 3 we give the corresponding result for longest paths and some corollaries.

\section{Proof of Theorem 8}
 We will need the following lemmas.

{\bf Lemma 9.} (Chv\'{a}tal [3]) {\it Let $G$ be a graph with degree sequence $d_1\le d_2\le\cdots\le d_n$ where $n\ge 3.$ If there is no integer $i$ with $1\le i<n/2$ such that $d_i\le i$ and $d_{n-i}<n-i,$ then $G$ is hamiltonian.}

Lemma 9 can also be found in [2, p.488].

{\bf Lemma 10.} (Dirac [4]) {\it Every $2$-connected $n$-vertex graph $G$ has a cycle of length at least ${\rm min}\{n, 2\delta(G)\}$.}

Let $S,T$ be disjoint subsets of $V(G)$, denote by $(S,T)$-path a path in $G$ from $S$ to $T$ has only endpoints in $S\cup T.$ If $S=\{x\}$ and $T=\{y\}$, let $P(x,y)$ denote an $(S,T)$-path and $|P(x,y)|$ denote the order of $P(x,y)$. Let $P$ be a path. If $u,v\in V(P)$, the closed interval $P[u, v]$ is defined to be the subpath of $P$ with endpoints $u$ and $v.$

{\bf Lemma 11.} {\it Let $G$ be a $k$-connected graph, and let $S,T$ be disjoint subsets of $V(G)$ with cardinality at least $k$. Then $G$ has $k$ pairwise disjoint $(S,T)$-paths.}

Lemma 11 can  be found in [15, p.174]. To prove our main result, we need the following lemma which is stated in [10] without proof.
This important lemma has been used in several papers. For completeness, we give a proof.

{\bf Lemma 12.} (Kopylov [10]) {\it Let $P$ be an $(x,y)$-path of length $m$ in a $2$-connected graph $G$. Then
$$
c(G)\geq {\rm min}\{m+1,\, d_P(x)+d_P(y)\}.
$$}

{\bf Proof.}
  Let $P=x_0x_1,\cdots x_m$ with $x_0=x$ and $x_m=y.$ Let $I=\{i|x_i\leftrightarrow x \}$ and  $J=\{j|x_j\leftrightarrow y \}.$ We distinguish three cases.

{\bf Case $1.$}\ There exists a vertex a vertex $x_i$, where $1\leq i\leq m-1$ such that $x_0$ is connected to $x_i$ and $x_m$ to $x_{i-1},$ then $xx_i\cup P[x_i,y]\cup yx_{i-1}\cup P[x_{i-1},x]$ is a cycle of length at least $m+1.$

{\bf Case $2.$}\ ${\rm Max}(I)> {\rm min}(J)$ and for any $i\in I, i-1\notin J.$ Chose $i\in I$ and $j\in J$ such that $i>j$ and $i-j$ is minimum. Denote by $k$ the order of the cycle $C=xx_i\cup P[x_i,y]\cup yx_j\cup P[x_j,x]$. We claim that $k\geq d_P(x)+d_P(y).$ First note that $N_P(x)\cup N_P(y)\cup \{x,y\}\subseteq V(C).$ Thus $k\geq d_p(y)+1.$ If $s\in I$ and $s\neq i,$ then $x_{s-1}\in V(C)$ and $s-1\notin J.$ It follows that $k\geq (d_P(y)+1)+(d_P(x)-1)=d_P(x)+d_P(y).$

{\bf Case $3.$}\ ${\rm Max}(I)\leq {\rm min}(J).$ Let $i={\rm max}(I)$ and $j={\rm min}(J).$ Denote $S=\{x,x_1,\cdots,x_i\}$ and $T=\{x_j,x_{j+1},\cdots,y\}.$ Since $G$ is $2$-connected, by Lemma 11, there exist two vertex disjoint $(S,T)$-paths $Q_1$ and $Q_2.$ Let $V(Q_1)\cap S=\{x_a\},$ $V(Q_1)\cap T=\{x_b\},$ $V(Q_2)\cap S=\{x_c\}$ and $V(Q_2)\cap T=\{x_d\}.$  Define $a',b',c',d'$ as follows, $a'$ be the smallest integer such that $a'\in I$ and $a'>a,$ $b'$ be the largest integer such that $b'\in J$ and $b'<b,$ $c'$ be the smallest integer such that $c'\in I$ and $c'>c$ and $d'$ be the largest integer such that $d'\in J$ and $d'<d.$ Denote $W=P[x_i,x_j].$

If one of $Q_1$ and $Q_2,$ say $Q_1,$ is vertex-disjoint from $W,$ the cycle
$$
xx_{a'}\cup P[x_{a'},x_{b'}]\cup x_{b'}y\cup P[y,x_b]\cup Q_1\cup P[x_a,x]
$$
contains $N_P(x)\cup N_P(y)\cup\{x,y\}.$ Hence $c(G)\geq d_P(x)+d_P(y)$.

Now suppose both $Q_1$ and $Q_2$ intersect $W.$ Let $x_q,x_r\in [V(Q_1)\cup V(Q_2)]\cap V(W)$ be such that the subscript $q$ is the smallest and $r$ is the largest. Without lost of generality, suppose $x_q\in V(Q_1).$

If $x_r\in Q_1,$ we have the cycle
$$
xx_{c'}\cup P[x_{c'},x_q]\cup Q_1[x_q,x_r] \cup P[x_r,x_{d'}]\cup x_{d'}y \cup P[y,x_d]\cup Q_2\cup P[x_c,x].
$$
This cycle contains $N_P(x)\cup N_P(y)\cup \{x,y\}.$ Hence $c(G)\geq d_P(x)+d_P(y)$.

If $x_r\in Q_2,$ we have the cycle
$$
xx_{c'}\cup P[x_{c'},x_q]\cup Q_1[x_q,x_b] \cup P[x_b,y]\cup yx_{b'} \cup P[x_{b'},x_r]\cup Q_2[x_r,x_c]\cup P[x_c,x].
$$
This cycle contains $N_P(x)\cup N_P(y)\cup \{x,y\}.$ Hence $c(G)\geq d_P(x)+d_P(y).$ All the possible cases have been considered, and hence the proof is complete. \hfill $\Box$

The following definition of {\it$t$-disintegration} of a graph is due to Kopylov [10].

{\bf Definition 1.} ($t$-disintegration of a graph, Kopylov [10]). Let $G$ be a graph and $t$ be a natural number. Delete all vertices of degree at most $t$ from $G$; for the resulting graph $G'$, we again delete all vertices of degree at most $t$ from $G'$. Iterating this process until we finally obtain a graph, denoted by $D(G; t)$, such that either $D(G; t)$ is a null graph or $\delta(D(G; t))\geq t+1.$ The graph $D(G; t)$ is called the {\it $(t+1)$-core of $G.$}

{\bf Definition 2.} Let $W$ be a set of vertices in a graph $G.$ $G$ is called {\it edge-maximal} with respect to the circumference and $W$ if for any $e\in E(\overline{G})$ such that $e$ has no endpoint in $W$, $c(G+e)> c(G).$

Now we are ready to prove Theorem 8.  We will use ideas from [10] (proof of Theorem 3), [11] (proof of Theorem 1.4) and
[13] (proof of Theorem 3.4). We also need to treat new situations, since more precise conditions are given in our problem.

{\bf Proof of Theorem 8.} It is easy to verify that the graphs $F_s(n,c,k)$ and $G_s(n,c,k)$ stated in Notation 1 are graphs of order $n$, circumference $c$ and minimum degree $k.$ The number of $K_s$ in $F_s(n,c,k)$ or $G_s(n,c,k)$ is $\varphi_s(n,c,k).$ Thus it remains to show that $\varphi_s(n,c,k)$ is an upper bound.

Let $Q$ be a graph of order $n$, circumference $c$ and minimum degree $k$. Choose a vertex $w$ of $Q$ with minimum degree. Suppose $Q$ is edge-maximal with respect to the circumference $c$ and $\{w\}.$ Thus each pair of non-adjacent vertices in $Q$ that does not contain $w$ is connected by a path of length at least $c.$ Let $t=\lfloor c/2\rfloor$ and $D(Q;t)$ denote the $(t+1)$-core of $Q$, i.e., the resulting graph of applying $t$-disintegration to $Q.$  For convenience, let $D=D(Q;t)$, we distinguish two cases.

{\bf Case 1.} $D$ is a null graph. In the $t$-disintegration process, denote $Q_0=Q$ and $Q_{i+1}=Q_i-x_i,0\leq i\leq n-1$ where $x_i$ is a vertex of degree at most $t$ in $Q_i$. Since $\delta(Q)\leq \lfloor\frac{c}{2}\rfloor=t$ by Lemma 10, We can always let $x_0=w.$  By the definition of $t$-disintegration, we have $d_{Q_i}(x_i)\leq t,1\leq i\leq n-t-1.$ Thus
$$
  N_s(Q)\leq\binom{k}{s-1}+(n-t-1)\binom{t}{s-1}+\binom{t}{s}\le g_s(n,c,k).
$$

{\bf Case 2.} $D$ is not a null graph. Let $d=|D|$, we claim that $D$ is a complete graph and $\delta (Q)=k \leq c-d+1$.

If there exist two vertices that are not adjacent in $D$, then in $Q$, there is a path of length at least $c$ with these vertices as its endpoints.
Among all nonadjacent pairs of vertices in $D$, choose $x, y\in V(D)$ such that $|P(x,y)|={\rm max}\{|P(u,v)|:u,v\in V(D), u\nleftrightarrow v\}.$
Note that $w$ is not in $V(D)$ and therefore $x,y\ne w.$ Let $P_1=P(x,y)$. We next show $N_D(x)\subseteq V(P_1)$ and $N_D(y)\subseteq V(P_1).$ If $x$ has a neighbor $z\in V(D)$ and $z\notin V(P_1)$ , then either $yz\in E(Q)$ and $zx\cup P_1\cup yz$ is a cycle of length at least $c+1$, or $yz\notin E(Q)$ and so $P_1\cup xz$ is a longer path. This contradicts the maximality of $P_1.$ Similar for $y$, we have $N_D(y)\subseteq V(P_1).$ Hence, by Lemma 12, $Q$ has a cycle of length at least ${\rm min}\{c+1,d_{P_1}(x)+d_{P_1}(y)\}\geq {\rm min}\{c+1,2(t+1)\}= c+1,$ a contradiction. Thus $D$ is a complete graph.

Suppose $k\geq c-d+2$, then $d\geq c-k+2.$ By the definition of $t$-disintegration, the minimum degree of $D$ is at least $t+1$, so we have $d\geq t+2.$ If $u\in V(Q)\backslash V(D)$, then $u$ is not adjacent to at least one vertex in $D.$  Choose $x\in V(Q)\backslash V(D)$ and $y\in V(D)$ such that $|P(x,y)|={\rm max}\{|P(u,v)|:u\in V(Q)\backslash V(D), v\in V(D), u\nleftrightarrow v\}.$ Denote $P_2=P(x,y)$. Now we will show that $P_2$ has length at least $c.$ Note that $w\in V(Q)\setminus V(D).$ We distinguish two cases. If $V(Q)\setminus V(D)=\{w\},$ we have $x=w$ and $|D|=n-1.$ We obtain a hamiltonian-path between $x$ and $y$ since $D$ is complete. Otherwise there exists an $(u,v)$-path with length at least $c.$ By the maximality of $P_2,$ we have the size of $P_2$ is at least $c.$ We claim that $N_Q(x)\subseteq V(P_2)$ and $N_D(y)\subseteq V(P_2).$  If $x$ has a neighbor $z\in V(D)$ and $z\notin V(P_2).$ Since $D$ is a complete graph, $zx\cup P_2\cup yz$ is a cycle of length at least $c+1$, a contradiction. If $x$ has a neighbor $z\in V(Q)\backslash V(D)$ and $z\notin V(P_2),$ then either $yz\in E(Q)$ and $zx\cup P_2\cup yz$ is a cycle of length at least $c+1$, or $yz\notin E(Q)$ and so $P_2\cup xz$ is a longer path. This contradicts the maximality of $P_2.$ Similar for $y$, we have $N_D(y)\subseteq V(P_2).$ Hence, by Lemma 12, there is a cycle with length at least ${\rm min}\{c+1,d_{P_2}(x)+d_{P_2}(y)\}\geq {\rm min}\{c+1,k+d-1\}\geq {\rm min}\{c+1,k+c-k+1\}= c+1,$ where the second inequality follows from $d\geq c-k+2,$ a contradiction. Thus $k\leq c-d+1.$

Apply $(c-d+1)$-disintegration to $Q$, and let $D'=D(Q;c-d+1)$ be the resulting graph. Recall that $d\geq t+2$. We have $k\leq c-d+1\leq c-t-1 \leq t$.  Then $w\notin V(D')$ and $D\subseteq D'.$ There are two cases.

(a) If $D'=D$, then $|D'|=|D|=d.$ By the definition of $(c-d+1)$-disintegration, we have
\begin{align*}
  N_s(Q)\leq \binom{k}{s-1}+(n-d-1)\binom{c-d+1}{s-1}+\binom{d}{s}&=\binom{k}{s-1}+\lambda_s(n,c,c-d+1)\\
  &\leq {\rm max}\{f_s(n,c,k),g_s(n,c,k)\},
\end{align*}
where $\lambda_s(n,c,x)=(n-c-2+x)\binom{x}{s-1}+\binom{c+1-x}{s}.$  The third inequality follows from the condition $k\leq c-d+1 \leq t$ and that the function $\lambda_s(n,c,x)$ is convex for $x\in[k,t].$

(b) Otherwise $D\neq D'.$ If $u\in V(D')\backslash V(D)$, then $u$ is not adjacent to at least one vertex in $D.$  Among all these nonadjacent pairs of vertices, choose $x\in V(D')\backslash V(D)$ and $y\in V(D)$ such that $|P(x,y)|={\rm max}\{|P(u,v)|: u\in V(D')\backslash V(D), v\in V(D), u\nleftrightarrow v\}.$ Denote $P_3=P(x,y)$ for convenience. As before, we have  $N_{D'}(x)\subseteq V(P_3)$ and $N_D(y)\subseteq V(P_3).$  By Lemma 12, there is a cycle with length at least ${\rm min}\{c+1, d_{P_3}(x)+d_{P_3}(y)\} \geq {\rm min}\{c+1,(c-d+2)+(d-1)\}=c+1$, where the second inequality follows from the condition $D'$ is the $(c-d+2)$-core of $Q$ and $|D|=d,$ a contradiction.

This completes the proof. \hfill $\Box$

\section{ The size version and further results}

Obviously, Theorem 2, Theorem 6 and Theorem 7 can be deduced from Theorem 8. By Theorem 8, we also have the following result.

{\bf Corollary 13.} {\it Let $\varphi(n,c,k)$ denote the maximum size of a $2$-connected nonhamiltonian graph of order $n$ with circumference $c$ and minimum degree $k.$  Then
$$
\varphi(n,c,k)= {\rm max}\{f(n,c,k),\, g(n,c,k)\}.
$$}

The following corollary follows from corollary 13 and Lemma 9.

{\bf Corollary 14.} {\it Let $\phi(n,k)$ denote the maximum size of a nonhamiltonian 2-connected graph of order $n$ and minimum degree $k.$ Then
$$
\phi(n,k)=\begin{cases} \binom{n-k}{2}+k^2\quad {\rm if}\,\,\,n\,\,\,{\rm is}\,\,\,{\rm odd}\,\,\,{\rm and}\,\,\,n\ge 6k-5\,\,\,{\rm or}\,\,\,
n\,\,\,{\rm is}\,\,\,{\rm even}\,\,\,{\rm and}\,\,\,n\ge 6k-8,\\
\frac{3n^2-8n+5}{8}+k\quad {\rm if}\,\,\,n\,\,\,{\rm is}\,\,\,{\rm odd}\,\,\,{\rm and}\,\,\,2k+1\le n\le 6k-7,\\
\frac{3n^2-10n+16}{8}+k\quad {\rm if}\,\,\,n\,\,\,{\rm is}\,\,\,{\rm even}\,\,\,{\rm and}\,\,\,2k+2\le n\le 6k-10.
\end{cases}
$$
If $n=6k-5$ or $n=6k-8,$ then $\phi(n,k)$ is attained by a graph $G$ if and only if $G=K_k\vee(K_{n-2k}+\overline{K_k})$ or $G=G(n,n-1,k).$
If $n$ is odd and $n\ge 6k-3$ or $n$ is even and $n\ge 6k-6,$ then $\phi(n,k)$ is attained by a graph $G$ if and only if $G=K_k\vee(K_{n-2k}+\overline{K_k}).$
If $n$ is odd and $2k+1\le n\le 6k-7$ or $n$ is even and $2k+2\le n\le 6k-10,$ then $\phi(n,k)$ is attained by a graph $G$ if and only if $G=G(n,n-1,k).$}

{\bf Proof.}
Let $\varphi(n,c,k)$ be defined as in corollary 13. It is easy to verify that
$$
{\rm max}\{\varphi(n,c,k):\,\, 2k\le c\le n-1\}=\phi(n,k).
$$
Now we determine the extremal graphs. Since the proof when $n$ is even is similar to the proof when $n$ is odd, we give only the proof of the latter.

Let $Q$ be a nonhamiltonian graph of order $n$ and minimum degree $k.$ The degree sequence of $Q$ is $d_1\le d_2\le\cdots\le d_n$ where $n\ge 3.$  By Lemma 9, there exists $i$ with $i\le(n-1)/2$ such that $d_i\le i$ and $d_{n-i}\le n-i-1.$ Assume that each vertex degree attains the value of the upper bound. We have degree sequence
$$
k,\, \underbrace{i,\ldots,i}_{i-1},\, \underbrace{n-i-1,\ldots,n-i-1}_{n-2i},\, \underbrace{n-2,\ldots,n-2}_{i-k},\, \underbrace{n-1,\ldots,n-1}_{k}.
$$
This degree sequence is graphical and the sum of all vertex degree is $n^2-(2i+1)n+3i^2-i+2k.$

If $n\ge 6k-5,$ suppose that $Q$ has size $\binom{n-k}{2}+k^2.$  We have
$$
e(Q)=\binom{n-k}{2}+k^2\le [n^2-(2i+1)n+3i^2-i+2k]/2, \eqno (1)
$$
where the inequality is equivalent to $(i-k)(2n-3i-3k+1)\le 0.$ Since $i\ge d_i\ge \delta(Q)=k,$ we obtain $i=k$ or $n\le (3i+3k-1)/2.$

If $i=k,$ equality holds in (1) and hence the degree sequence of $Q$ is
$$
\underbrace{k,\ldots,k}_{k},\, \underbrace{n-k-1,\ldots,n-k-1}_{n-2k},\, \underbrace{n-1,\ldots,n-1}_{k},
$$
implying that $Q=K_k\vee(K_{n-2k}+\overline{K_k}).$

Now suppose $i\neq k.$ Then we have $n\le (3i+3k-1)/2.$ If $i\le (n-3)/2,$ then $n\le 6k-11,$ contradicting our assumption $n\ge 6k-5.$ Thus $i=(n-1)/2.$
We have $n\leq 6k-5.$ Since $n\ge 6k-5,$ we have $n=6k-5.$  Hence the degree sequence of $Q$ is
$$
k,\, \underbrace{\frac{n-1}{2},\ldots,\frac{n-1}{2}}_{\frac{n-1}{2}},\, \underbrace{n-2,\ldots,n-2}_{\frac{n-1}{2}-k},\, \underbrace{n-1,\ldots,n-1}_{k},
$$
implying that $Q=G(n,n-1,k).$

Otherwise $2k+1\le n\le 6k-7.$ Suppose that $Q$ has size $\frac{3n^2-8n+5}{8}+k.$ We have
$$
e(Q)=\frac{3n^2-8n+5}{8}+k\le [n^2-(2i+1)n+3i^2-i+2k]/2, \eqno (2)
$$
where the inequality is equivalent to $0\le (6i-n-5)(2i-n+1).$ It follows that $i\le (n+5)/6$ or $i\ge (n-1)/2$. Since $(n-1)/2\ge i\ge d_i\ge \delta(Q)=k$ and $n\le 6k-7,$  we obtain $i=(n-1)/2.$ Hence the degree sequence of $Q$ is
$$
k,\, \underbrace{\frac{n-1}{2},\ldots,\frac{n-1}{2}}_{\frac{n-1}{2}},\, \underbrace{n-2,\ldots,n-2}_{\frac{n-1}{2}-k},\, \underbrace{n-1,\ldots,n-1}_{k}.
$$
Let $u_i\, (i=1,\ldots, k),$ $v_j\, (j=1,\ldots, (n-1)/2-k),$ $w_s\, (s=1,\ldots,(n-1)/2)$ and $x$ be the vertices of $Q$ where each $u_i$ has degree $n-1,$
each $v_j$ has degree $n-2,$ each $w_s$ has degree $(n-1)/2$ and $x$ has degree $k.$ Then each $u_i$ is adjacent to every other vertex. The vertex $x$ has exactly the $u_i's$ as its neighbors. Let $H=Q-\{x,u_1,\ldots, u_k\}.$ Clearly each $v_j$ is a dominating vertex of $H.$ Now every $w_s$ has exactly the $v_j's$
as its neighbors and $w_1,\ldots, w_{(n-1)/2}$ form an independent set. This shows that that $Q=G(n,n-1,k).$ This completes the proof. \hfill $\Box$

A {\it detour} of a graph $G$ is a longest path in $G.$ The {\it detour order} of $G$ is the number of vertices in a detour of $G.$ The following trick is well-known (e.g. [4, p.166] or [15, p.292]).

{\bf Lemma 15.} {\it Let $G$ be a graph and denote $H=G\vee K_1.$ Then the detour order of $G$ is $p$ if and only if the circumference of $H$ is $p+1.$}

The following corollary follows from corollary 13 and Lemma 15 immediately. A graph is called {\it traceable} if it has a Hamilton path; otherwise it is
{\it nontraceable.}

{\bf Corollary 16.} {\it Let $\psi(n,p,k)$ denote the maximum size of a nontraceable graph of order $n,$ detour order $p$ and minimum degree $k.$  Then
$$
\psi(n,p,k)= {\rm max}\{f(n,p-1,k),\, g(n,p-1,k)\}.
$$ }

For positive integers $n,\, p$ and $k,$ let $t=\lfloor (p+1)/2\rfloor.$  $G'(n,p,k)$ denotes the graph obtained from
$K_{t-1}\vee (K_{p+2-2t}+\overline{K_{n-p+t-1}})$ by deleting $t-1-k$ edges that are incident to one common vertex in $\overline{K_{n-p+t-1}}.$
 Let $F'(n,p,k)=K_k\vee (K_{p-2k}+\overline{K_{n-p+k}}).$ Both $F'(n,p,k)$ and $G'(n,p,k)$ have order $n$, minimum degree $k$ and detour order $p$. Since $F'(n,p,k)$ and $G'(n,p,k)$ have sizes $f(n,p-1,k)$ and $g(n,p-1,k)$ respectively, the bound in corollary 16 is best possible.

{\bf Notation 3.} Fix $n-1 \geq c\geq 2k\ge 4, t=\lfloor c/2\rfloor.$  For $K_t\vee(K_{c+1-2t}+\overline{K_{n-c-1+t}}),$  let $W$ be the vertex set with $q$ vertices in $\overline{K_{n-c-1+t}}.$ For every vertex in $W,$ delete $t-k$ edges that are incident to it. Let $G(n,c,k,q)$ denote the resulting graph. Denote by $g_s(n,c,k,q)$ the number of $s$-cliques in $G(n,c,k,q)$; more precisely,
$$
  g_s(n,c,k,q)=q\binom{k}{s-1}+(n-q-c-1+t)\binom{t}{s-1}+\binom{c+1-t}{s}.
$$
Note that in the proof of Theorem 8, if we replace the vertex $w$ with vertices set $W$ that has $q$ minimum degree vertices, where $q\leq k,$ we will have the following result.

{\bf Corollary 17.}  {\it Let $G$ be a $2$-connected nonhamiltonian graph of order $n$ with circumference $c$ and at least $q$ vertices with minimum degree $k,$ $q\leq k.$ Then
$$
N_s(G)\leq {\rm max}\{f_s(n,c,k),g_s(n,c,k,q)\}.
$$ }

\vskip 5mm
{\bf Acknowledgement.} The author is grateful to Professor Xingzhi Zhan for suggesting the problems and for helpful discussions, and to one referee for her (his) kind suggestions. This research  was partially supported by the NSFC grants 11671148 and Science and Technology Commission of Shanghai Municipality (STCSM) grant 18dz2271000.


\begin{thebibliography}{99}
\bibitem{1} J.A. Bondy, Variations on the hamiltonian theme, Canad. Math. Bull., 15(1972), no. 1, 57-62.
\bibitem{3} J.A. Bondy and U.S.R. Murty, Graph Theory, GTM 244, Springer, 2008.
\bibitem{4} V. Chv\'{a}tal, On Hamilton's ideals, J. Combin. Theory, Ser. B, 12(1972), 163-168.
\bibitem{2} G.A. Dirac, Some theorems on abstract graph, proc. Lond. Math. Soc., 2(1952), 69-81.
\bibitem{6} P. Erd\H{o}s, Remarks on a paper of P\'{o}sa, Magyar Tud. Akad. Math. Kutat\'{o}. Int. K\"{o}zl., 7(1962), 227-229.
\bibitem{11} P. Erd\H{o}s, T. Gallai, On maximal paths and circuits of graphs, Acta Math. Acad. Sci. Hung., 10(1959), 337-356.
\bibitem{14}G. Fan, X. Lv and P. Wang, Cycles in 2-connected graphs, J. Combin. Theory Ser. B, 92(2004), 379-394.
\bibitem{5} Z. F\"{u}redi, A. Kostochka and R. Luo, Extensions of a theorem of Erd\H{o}s on nonhamiltonian graphs, J. Graph Theory, 89(2018), 176-193.
\bibitem{15} Z. F\"{u}redi, A. Kostochka and J. Verstra\"{e}te, Stability in the Erd\H{o}s-Gallai theorems on cycles and paths, J. Combin. Theory Ser. B, 121(2016), 197-228.
\bibitem{9} G.N. Kopylov, Maximal paths and cycles in a graph, Dokl. Akad. Nauk SSSR, 234(1977), 19-21. English translation: Soviet Math. Dokl. 18(1977), 593-596.
\bibitem{7} R. Luo, The maximum number of cliques in graphs without long cycles, J. Combin. Theory Ser. B, 128(2017), 219-226.
\bibitem{16} J. Ma, B. Ning, Stability results on the circumference of a graph, Combinatorica, 40(1)(2020), 105-147.
\bibitem{12} B. Ning, X. Peng, Extensions of the Erd\H{o}s-Gallai theoremand Luo theorem, Comb. Probab. Comput., 29(2020), 128-136.
\bibitem{8} O. Ore, Arc coverings of graphs, Ann. Mat. Pura Appl., 55(1961), 315-321.
\bibitem{10} D.B. West, Introduction to Graph Theory, Prentice Hall, Inc., 1996.
\bibitem{13} D.R. Woodall, Maximal circuits of graphs I, Acta Math. Acad. Sci. Hungar., 28(1976), 77-80.
\end{thebibliography}
\end{document}